\documentclass[12pt,american]{article}
\usepackage[T1]{fontenc}
\usepackage[latin9]{inputenc}
\usepackage{color}
\usepackage{babel}
\usepackage{amsmath}
\usepackage{amsthm}
\usepackage{amssymb}
\usepackage[unicode=true,
 bookmarks=true,bookmarksnumbered=true,bookmarksopen=false,
 breaklinks=true,pdfborder={0 0 1},backref=false,colorlinks=true]
 {hyperref}
\hypersetup{pdftitle={Cesaro Limits for Fractional Dynamics},
 pdfauthor={Jos{\'e} L.~da Silva and Yuri G. Kondratiev},
 pdfsubject={Cesaro Limits for Fractional Dynamics},
 pdfkeywords={Dynamical Systems, random time, subordinator, asymptotic behavior},
 linkcolor=blue,citecolor=black,urlcolor=black,filecolor=black}

\makeatletter
\theoremstyle{plain}
\newtheorem{thm}{\protect\theoremname}[section]
\theoremstyle{definition}
\newtheorem{example}[thm]{\protect\examplename}
\theoremstyle{remark}
\newtheorem{rem}[thm]{\protect\remarkname}
\theoremstyle{plain}
\newtheorem{lem}[thm]{\protect\lemmaname}

\usepackage{babel}

\usepackage{babel}
\usepackage{soul}
\theoremstyle{plain}
\usepackage{units}
\usepackage{soul}
\usepackage{textcomp}
\usepackage{mathrsfs}
\usepackage{bbm}
\usepackage{pdfsync}
\usepackage{datetime}\usepackage{currfile}

\numberwithin{equation}{section}

\newcommand{\R}{{\mathbb R}}
\newcommand{\X}{{\R^d}}

\newcommand{\E}{{\mathbb E}}

\newcommand{\om}{\omega}

\makeatother

\providecommand{\examplename}{Example}
\providecommand{\lemmaname}{Lemma}
\providecommand{\remarkname}{Remark}
\providecommand{\theoremname}{Theorem}

\begin{document}
\title{Cesaro Limits for Fractional Dynamics}
\author{\textbf{Yuri Kondratiev}\\
 Department of Mathematics, University of Bielefeld,\\
 D-33615 Bielefeld, Germany,\\
 Dragomanov University, Kiev, Ukraine\\
 Email: kondrat@math.uni-bielefeld.de\and\textbf{Jos{\'e} Lu{\'i}s
da Silva},\\
 CIMA, University of Madeira, Campus da Penteada,\\
 9020-105 Funchal, Portugal.\\
 Email: joses@staff.uma.pt}
\date{\today}
\maketitle
\begin{abstract}
We study the asymptotic behavior of random time changes of dynamical
systems. As random time changes we propose three classes which exhibits
different patterns of asymptotic decays. The subordination principle
may be applied to study the asymptotic behavior of the random time
dynamical systems. It turns out that for the special case of stable
subordinators explicit expressions for the subordination are known
and its asymptotic behavior are derived. For more general classes
of random time changes explicit calculations are essentially more
complicated and we reduce our study to the asymptotic behavior of
the corresponding Cesaro limit.


\end{abstract}
\tableofcontents{}

\section{Introduction}

In this paper we will deal with Markov processes or dynamical systems
in $\X$. These processes or dynamics starting from $x\in\X$, denote
by $X^{x}(t)$, $t\geq0$, have associated evolution equations on
$\X$. In the Markov case we define for suitable $f:\X\longrightarrow\R$
the function $u(t,x)=\E[f(X^{x}(t))]$ which satisfied the Kolmogorov
equation 
\[
\frac{\partial}{\partial t}u(t,x)=Lu(t,x),
\]
where $L$ is the generator of the \textcompwordmark Markov process.

For a dynamical system we introduce $u(t,x)=f(X^{x}(t))$. Then this
function is the solution of the Liouville equation 
\[
\frac{\partial}{\partial t}u(t,x)=Lu(t,x),
\]
where now $L$ is the Liouville operator for the dynamical system,
see e.g. \cite{KS2020a}.

Let $S(t)$, $t\ge0$ be a subordinator and $E(t)$, $t\ge0$ denotes
the inverse subordinator, that is, for each $t\ge0$, $E(t):=\inf\{s>0\mid S(s)>t\}$.
This random process we consider as a random time and assume to be
independent of $X^{x}(t)$. Define a random process $Y^{x}$ by\textcompwordmark
\[
Y^{x}(t,\om):=X^{x}(E(t,\om)).
\]
Then as above we may introduce 
\[
u^{E}(t,x)=\E[f(Y^{x}(t))].
\]
For both Markov and dynamical system cases this function satisfies
the evolution equations 
\[
D_{t}^{E}u^{E}(t,x)=Lu^{E}(t,x)
\]
where $L$ is the Kolmogorov or Liouville operator correspondingly.
Here $D_{t}^{E}$ is a generalized fractional time derivative corresponding
to the inverse subordinator $E(t)$, see Section~\ref{sec:RTP} below
for details, in particular the definition in \eqref{eq:general-derivative}.
The main relation which is true for both cases is the following subordination
formula:

\begin{equation}
u^{E}(t,x)=\int_{0}^{\infty}u(\tau,x)G_{t}(\tau)d\tau,\label{subor}
\end{equation}
where $G_{t}(\tau)$ is the density of the inverse subordinator $E(t)$,
see, e.g., \cite{Toaldo2015}, \cite{KS2020a} and especially the
book \cite{Meerschaert2012}. This formula which relates the solutions
of the evolution equations with usual and fractional derivatives plays
an important role in the study of dynamics with random times. Note
that there exist such relations between random times, fractional equations
and subordination in the framework of physical models, see, e.g.,
\cite{Mura_Taqqu_Mainardi_08}.

The goal of this paper is to study and analyze the asymptotic behavior
of two elementary dynamical system after the random time change, namely
$u(t,x)=e^{-at}$, $a>0$ and $u(t,x)=t^{n}$, $n\ge0$. Here the
dynamical system are considered as a deterministic Markov processes.
For particular classes of random times the subordination formula \eqref{subor}
is evaluated explicitly. This is true, for example, in the case of
inverse stable subordinators. For a general inverse subordinator the
properties of the density $G_{t}(\tau)$ are unknown and the evaluation
of \eqref{subor} is not possible. Actually, it is a long standing
open problem in the theory of stochastic processes.

We propose an alternative approach to study the asymptotic behavior
of $u^{E}(t,x)$. More precisely, we consider Cesaro limits (the asymptotic
of the Cesaro mean of $u^{E}(t,x)$, see \eqref{Cesaro-mean} below)
of $u^{E}(t,x)$ using the subordination formula representation \eqref{subor}
together with the Fe< ller--Karamata Tauberian theorem, see Theorem~\ref{thm:FKT-LST}.
For many classes of random times this approach leads to a precise
asymptotic behavior. In this paper we investigate three classes of
random time change, denote by \ref{eq:C1}, \ref{eq:C2}, and \ref{eq:C3},
see Section~\ref{sec:RTP}, which exhibits different patterns of
decays of the Cesaro limit of $u^{E}(t,x)$. We would like to emphasize
that for particular classes of random times, namely inverse stable
subordinators, the asymptotic of $u^{E}(t,x)$ which may be computed
explicitly, coincides with the Cesaro limit. For other classes of
random times the Cesaro limit gives one possible characteristic of
the asymptotic for $u^{E}(t,x)$. To the best of our knowledge at
the present time no other information on the asymptotic of $u^{E}(t,x)$
is known for a general subordinator.

The remaining of the paper is organized as follows. In Section \ref{sec:RTP}
we introduce three classes (\ref{eq:C1}, \ref{eq:C2}, and \ref{eq:C3})
of subordinator processes which serves as random times. These classes
are given in terms of their local behavior of the Laplace exponent
at $\lambda=0$. In addition, we state the main results of the paper.
Section \ref{sec:ISS} is a preparation for the more general study
of the asymptotic of the subordination in Section \ref{sec:General-Classes}.
More precisely, we investigate in detail the special case of the inverse
stable subordinator where explicit expressions are known. Hence, the
expression for the subordination \eqref{subor} is derived (for the
two dynamical systems considered above) as well as their Cesaro limit.
It turns out that both asymptotic for $u^{E}(t,x)$ (the explicit
calculations and Cesaro limit) are the same. Finally in Section \ref{sec:General-Classes}
we study the Cesaro limit for the general classes \ref{eq:C1}, \ref{eq:C2},
and \ref{eq:C3} of random time changes.

\section{Random Times Processes}

\label{sec:RTP}In this section we introduce three classes of subordinators
which serves as random times processes. More precisely, the random
times corresponds to the inverse of subordinator processes whose Laplace
exponent satisfies certain conditions, see below for details. The
simplest example in class \ref{eq:C1} below, is the well known $\alpha$-stable
subordinators whose inverse processes are well studied in the literature,
see for example \cite{Bingham1971} or \cite{Feller71}.

The classes of processes to be introduced which serve as random times
have a connection with the concept of general fractional derivatives
(see \cite{Kochubei11} for details and applications to fractional
differential equations) associated to an admissible kernels $k\in L_{\mathrm{loc}}^{1}(\mathbb{R}_{+})$
which is characterized in terms of their Laplace transforms $\mathcal{K}(\lambda)$
as $\lambda\to0$, see assumption (H) below.

\subsection{Definitions and Main Assumptions }

Let $S=\{S(t),\;t\ge0\}$ be a subordinator without drift starting
at zero, that is, an increasing L{\'e}vy process starting at zero,
see \cite{Bertoin96} for more details. The Laplace transform of $S(t)$,
$t\ge0$ is expressed in terms of a Bernstein function $\Phi:[0,\infty)\longrightarrow[0,\infty)$
(also known as Laplace exponent) by 
\[
\mathbb{E}(e^{-\lambda S(t)})=e^{-t\Phi(\lambda)},\quad\lambda\ge0.
\]
The function $\Phi$ admits the L{\'e}vy-Khintchine representation
\begin{equation}
\Phi(\lambda)=\int_{(0,\infty)}(1-e^{-\lambda\tau})\,\mathrm{d}\sigma(\tau),\label{eq:Levy-Khintchine}
\end{equation}
where the measure $\sigma$ (called L{\'e}vy measure) has support
in $[0,\infty)$ and fulfills 
\begin{equation}
\int_{(0,\infty)}(1\wedge\tau)\,\mathrm{d}\sigma(\tau)<\infty.\label{eq:Levy-condition}
\end{equation}
In what follows we assume that the L{\'e}vy measure $\sigma$ satisfy
\begin{equation}
\sigma\big((0,\infty)\big)=\infty.\label{eq:Levy-massumption}
\end{equation}
Using the L{\'e}vy measure $\sigma$ we define the kernel $k$ as
follows 
\begin{equation}
k:(0,\infty)\longrightarrow(0,\infty),\;t\mapsto k(t):=\sigma\big((t,\infty)\big).\label{eq:k}
\end{equation}
Its Laplace transform is denoted by $\mathcal{K}$, that is, for any
$\lambda\ge0$ one has 
\begin{equation}
\mathcal{K}(\lambda):=\int_{0}^{\infty}e^{-\lambda t}k(t)\,\mathrm{d}t.\label{eq:LT-k}
\end{equation}
The relation between the function $\mathcal{K}$ and the Laplace exponent
$\Phi$ is given by 
\begin{equation}
\Phi(\lambda)=\lambda\mathcal{K}(\lambda),\quad\forall\lambda\ge0.\label{eq:Laplace-exponent}
\end{equation}

We make the following assumption on the Laplace exponent $\Phi(\lambda)$
of the subordinator $S$. 
\begin{description}
\item [{(H)}] $\Phi$ is a complete Bernstein function (that is, the L{\'e}vy
measure $\sigma$ is absolutely continuous with respect to the Lebesgue
measure) and the functions $\mathcal{K}$, $\Phi$ satisfy 
\begin{equation}
\mathcal{K}(\lambda)\to\infty,\text{ as \ensuremath{\lambda\to0}};\quad\mathcal{K}(\lambda)\to0,\text{ as \ensuremath{\lambda\to\infty}};\label{eq:H1}
\end{equation}
\begin{equation}
\Phi(\lambda)\to0,\text{ as \ensuremath{\lambda\to0}};\quad\Phi(\lambda)\to\infty,\text{ as \ensuremath{\lambda\to\infty}}.\label{eq:H2}
\end{equation}
\end{description}
\begin{example}[$\alpha$-stable subordinator]
\label{exa:alpha-stable1}A classical example of a subordinator $S$
is the so-called $\alpha$-stable process with index $\alpha\in(0,1)$.
Specifically, a subordinator is $\alpha$-stable if its Laplace exponent
is 
\[
\Phi(\lambda)=\lambda^{\alpha}=\int_{0}^{\infty}(1-e^{-\lambda\tau})\frac{\alpha\tau^{-1-\alpha}}{\Gamma(1-\alpha)}\,\mathrm{d}\tau.
\]
In this case it follows that the L{\'e}vy measure is $\mathrm{d}\sigma_{\alpha}(\tau)=\frac{\alpha}{\Gamma(1-\alpha)}\tau^{-(1+\alpha)}\,\mathrm{d}\tau$.
The corresponding kernel $k_{\alpha}$ has the form $k_{\alpha}(t)=g_{1-\alpha}(t):=\frac{t^{-\alpha}}{\Gamma(1-\alpha)}$,
$t\ge0$ and its Laplace transform is $\mathcal{K}_{\alpha}(\lambda)=\lambda^{\alpha-1}$,
$\lambda>0$. 
\end{example}

\begin{example}[Sum of two stable subordinators]
\label{exa:sum-two-stables}Let $0<\alpha<\beta<1$ be given and
$S_{\alpha,\beta}(t)$, $t\ge0$ the driftless subordinator with Laplace
exponent given by 
\[
\Phi_{\alpha,\beta}(\lambda)=\lambda^{\alpha}+\lambda^{\beta}.
\]
It is clear from Example \ref{exa:alpha-stable1} that the corresponding
L{\'e}vy measure $\sigma_{\alpha,\beta}$ is the sum of two L{\'e}vy
measures, that is, 
\[
\mathrm{d}\sigma_{\alpha,\beta}(\tau)=\mathrm{d}\sigma_{\alpha}(\tau)+\mathrm{d}\sigma_{\beta}(\tau)=\frac{\alpha}{\Gamma(1-\alpha)}\tau^{-(1+\alpha)}\,\mathrm{d}\tau+\frac{\beta}{\Gamma(1-\beta)}\tau^{-(1+\beta)}\,\mathrm{d}\tau.
\]
Then the associated kernel $k_{\alpha,\beta}$ is 
\[
k_{\alpha,\beta}(t):=g_{1-\alpha}(t)+g_{1-\beta}(t)=\frac{t^{-\alpha}}{\Gamma(1-\alpha)}+\frac{t^{-\beta}}{\Gamma(1-\beta)},\;t>0
\]
and its Laplace transform is $\mathcal{K}_{\alpha,\beta}(\lambda)=\mathcal{K}_{\alpha}(\lambda)+\mathcal{K}_{\beta}(\lambda)=\lambda^{\alpha-1}+\lambda^{\beta-1}$,
$\lambda>0$. 
\end{example}

Let $E$ be the inverse process of the subordinator $S$, that is,
\begin{equation}
E(t):=\inf\{s>0\mid S(s)>t\}=\sup\{s\ge0\mid S(s)\le t\}.\label{eq:inverse-sub}
\end{equation}
For any $t\ge0$ we denote by $G_{t}(\tau)$, $\tau\ge0$ the marginal
density of $E(t)$ or, equivalently 
\[
G_{t}(\tau)\,\mathrm{d}\tau=\frac{\partial}{\partial\tau}P(E(t)\le\tau)\,\mathrm{d}\tau=\frac{\partial}{\partial\tau}P(S(\tau)\ge t)\,\mathrm{d}\tau=-\frac{\partial}{\partial\tau}P(S(\tau)<t)\,\mathrm{d}\tau.
\]

The density $G_{t}(\tau)$ is the main object in our considerations
below. Therefore, in what follows, we collect the most important properties
of $G_{t}(\tau)$ needed in the next sections. 
\begin{rem}
\label{rem:distr-alphastab-E}If $S$ is the $\alpha$-stable process,
$\alpha\in(0,1)$, then the inverse process $E(t)$, has Laplace transform
(cf.~Prop.~1(a) in \cite{Bingham1971} or \cite{Feller71}) given
by 
\begin{equation}
\mathbb{E}(e^{-\lambda E(t)})=\int_{0}^{\infty}e^{-\lambda\tau}G_{t}(\tau)\,\mathrm{d}\tau=\sum_{n=0}^{\infty}\frac{(-\lambda t^{\alpha})^{n}}{\Gamma(n\alpha+1)}=E_{\alpha}(-\lambda t^{\alpha}),\label{eq:LT-ISS-density}
\end{equation}
where $E_{\alpha}$ is the Mittag-Leffler function. It follows from
the asymptotic behavior of the function $E_{\alpha}$ that $\mathbb{E}(e^{-\lambda E(t)})\sim Ct^{-\alpha}$
as $t\to\infty$. It is possible to find explicitly the density $G_{t}(\tau)$
in this case using the completely monotonic property of the Mittag-Leffler
function $E_{\alpha}$. It is given in terms of the Wright function
$W_{\mu,\nu}$, namely $G_{t}(\tau)=t^{-\alpha}W_{-\alpha,1-\alpha}(\tau t^{-\alpha})$,
see \cite{Gorenflo1999} for more details. 
\end{rem}

For a general subordinator, the following lemma determines the $t$-Laplace
transform of $G_{t}(\tau)$, with $k$ and $\mathcal{K}$ given in
\eqref{eq:k} and \eqref{eq:LT-k}, respectively. For the proof see
\cite{Kochubei11} or Proposition~3.2 in \cite{Toaldo2015}. 
\begin{lem}
\label{lem:t-LT-G}The $t$-Laplace transform of the density $G_{t}(\tau)$
is given by 
\begin{equation}
\int_{0}^{\infty}e^{-\lambda t}G_{t}(\tau)\,\mathrm{d}t=\mathcal{K}(\lambda)e^{-\tau\lambda\mathcal{K}(\lambda)}.\label{eq:LT-G-t}
\end{equation}
The double ($\tau,t$)-Laplace transform of $G_{t}(\tau)$ is 
\begin{equation}
\int_{0}^{\infty}\int_{0}^{\infty}e^{-p\tau}e^{-\lambda t}G_{t}(\tau)\,\mathrm{d}t\,\mathrm{d}\tau=\frac{\mathcal{K}(\lambda)}{\lambda\mathcal{K}(\lambda)+p}.\label{eq:double-Laplace}
\end{equation}
\end{lem}

Here we would like to make the connection of the above abstract framework
with general fractional derivatives. For any $\alpha\in(0,1)$ the
Caputo-Dzhrbashyan fractional derivative of order $\alpha$ of a function
$u$ is defined by (see e.g., \cite{KST2006} and references therein)
\begin{equation}
\big(\mathbb{D}_{t}^{\alpha}u\big)(t)=\frac{d}{dt}\int_{0}^{t}k_{\alpha}(t-\tau)u(\tau)\,\mathrm{d}\tau-k_{\alpha}(t)u(0),\quad t>0,\label{eq:Caputo-derivative}
\end{equation}
where $k_{\alpha}$ is given in Example \ref{exa:alpha-stable1},
that is, $k_{\alpha}(t)=g_{1-\alpha}(t)=\frac{t^{-\alpha}}{\Gamma(1-\alpha)}$,
$t>0$. In general, starting with a subordinator $S$ and the kernel
$k\in L_{\mathrm{loc}}^{1}(\mathbb{R}_{+})$ as given in \eqref{eq:k},
we may define a differential-convolution operator by 
\begin{equation}
\big(\mathbb{D}_{t}^{(k)}u\big)(t)=\frac{d}{dt}\int_{0}^{t}k(t-\tau)u(\tau)\,\mathrm{d}\tau-k(t)u(0),\;t>0.\label{eq:general-derivative}
\end{equation}
The operator $\mathbb{D}_{t}^{(k)}$ is also known as general fractional
derivative and its applications to convolution-type differential equations
was investigated in \cite{Kochubei11}. 
\begin{example}[Distributed order derivative]
\label{exa:distr-order-deriv}Consider the kernel $k$ defined by
\begin{equation}
k(t):=\int_{0}^{1}g_{\alpha}(t)\,\mathrm{d}\alpha=\int_{0}^{1}\frac{t^{\alpha-1}}{\Gamma(\alpha)}\,\mathrm{d}\alpha,\quad t>0.\label{eq:distributed-kernel}
\end{equation}
Then it is easy to see that 
\[
\mathcal{K}(\lambda)=\int_{0}^{\infty}e^{-\lambda t}k(t)\,\mathrm{d}t=\frac{\lambda-1}{\lambda\log(\lambda)},\quad\lambda>0.
\]
The corresponding differential-convolution operator $\mathbb{D}_{t}^{(k)}$
is called distributed order derivative, see \cite{Atanackovic2009,Daftardar-Gejji2008,Hanyga2007,Kochubei2008,Gorenflo2005,Meerschaert2006}
for more details and applications. 
\end{example}

We say that the functions $f$ and $g$ are \emph{asymptotically equivalent
at infinity}, and denote $f(x)\sim g(x)$ as $x\to\infty$, meaning
that 
\[
\lim_{x\to\infty}\frac{f(x)}{g(x)}=1.
\]
We say that a function $L$ is \emph{slowly varying at infinity} (see
\cite{Feller71,Seneta1976}) if 
\[
\lim_{x\to\infty}\frac{L(\lambda x)}{L(x)}=1,\quad\mathrm{for\;any}\;\lambda>0.
\]
Below $C$ is constant whose value is unimportant and may change from
line to line.

In the following we consider three classes of admissible kernels $k\in L_{\mathrm{loc}}^{1}(\mathbb{R}_{+})$,
characterized in terms of their Laplace transforms $\mathcal{K}(\lambda)$
as $\lambda\to0$ (i.e., as local conditions): 
\begin{equation}
\mathcal{K}(\lambda)\sim\lambda^{\alpha-1},\quad0<\alpha<1.\tag*{(C1)}\label{eq:C1}
\end{equation}
\begin{equation}
\mathcal{K}(\lambda)\sim\lambda^{-1}L\left(\frac{1}{\lambda}\right),\quad L(y):=C\log(y)^{-1},\;C>0.\tag*{(C2)}\label{eq:C2}
\end{equation}
\begin{equation}
\mathcal{K}(\lambda)\sim\lambda^{-1}L\left(\frac{1}{\lambda}\right),\quad L(y):=C\log(y)^{-1-s},\;s>0,\;C>0.\tag*{(C3)}\label{eq:C3}
\end{equation}
We would like to emphasize that these three classes of kernels leads
to different type of differential-convolution operators. In particular,
the Caputo-Djrbashian fractional derivative \ref{eq:C1} and distributed
order derivatives \ref{eq:C2}, \ref{eq:C3}. Moreover, it is simple
to check that the class of subordinators from Example \ref{exa:sum-two-stables}
falls into the class \ref{eq:C1} above.
\begin{rem}
\label{rem:two-estimates}The asymptotic behavior of the function
$f(t)$ as $t\to\infty$ may be determined, under certain conditions,
by studying the behavior of its Laplace transform $\tilde{f}(\lambda)$
as $\lambda\to0$, and vice versa. An important situation where such
a correspondence holds is described by the Feller--Karamata Tauberian
(FKT) theorem. 
\end{rem}

We state below a version of the FKT theorem which suffices for our
purposes, see the monographs \cite[Sec.~1.7]{Bingham1987} and \cite[XIII, Sec.~1.5]{Feller71}
for a more general version and proofs. 
\begin{thm}[Feller--Karamata Tauberian]
\label{thm:FKT-LST}Let $U:[0,\infty)\longrightarrow\mathbb{R}$
be a monotone non-decreasing right-continuous function such that 
\[
w(\lambda):=\int_{0}^{\infty}e^{-\lambda t}\,\mathrm{d}U(t)<\infty,\quad\forall\lambda>0.
\]
If $L$ is a slowly varying function and $C,\rho\ge0$, then the following
are equivalent 
\begin{equation}
U(t)\sim\frac{C}{\Gamma(\rho+1)}t^{\rho}L(t)\quad\mathrm{as}\;t\to\infty,\label{eq:asymp-U}
\end{equation}
\begin{equation}
w(\lambda)\sim C\lambda^{-\rho}L\left(\frac{1}{\lambda}\right)\quad\mathrm{as}\;\lambda\to0^{+}.\label{eq:asym-w}
\end{equation}
When $C=0$, (\ref{eq:asymp-U}) is to be interpreted as $U(t)=o(t^{\rho}L(t))$;
similarly for (\ref{eq:asym-w}). 
\end{thm}

\subsection{Statement of the Main Results}

In Section \ref{sec:ISS} and \ref{sec:General-Classes} we will focus
our attention on deriving the asymptotic behavior of the subordination
$u^{E}(t,x)$ given in \eqref{subor} for the inverse stable subordinator
as well as for the classes \ref{eq:C1}, \ref{eq:C2}, and \ref{eq:C3}
given above. On one hand, the results concerning the inverse stable
subordinator as a random time are well understood, due to the fact
that the Laplace transform (in $\tau)$ of the density $G_{t}(\tau)$
is known (cf.~Remark \ref{rem:distr-alphastab-E}). On the other
hand, for a general subordinator much less information about $G_{t}(\tau)$
is known and explicit results for the subordination $u^{E}(t,x)$
are not available. In order to get around this problem, and motivated
by the results of Section \ref{sec:ISS}, we study the Cesaro limit
of $u^{E}(t,x)$ for the general classes of random times.

With the above considerations we are ready to state our main results. 
\begin{thm}
\label{thm:ISS}Let $u^{E}(t,x)$ be the subordination by the density
$G_{t}(\tau)$ associated to the inverse stable subordinator. Denote
by $M_{t}(u^{E}(\cdot,x)):=\frac{1}{t}\int_{0}^{t}u^{E}(s,x)\,\mathrm{d}s$
the Cesaro mean of $u^{E}(t,x)$. 

\begin{enumerate}
\item If $u(t,x)=t^{n}$, $n\ge0$, then the asymptotic behavior of $u^{E}(t,x)$
coincides with the Cesaro limit and is equal to 
\[
Ct^{n\alpha}\quad\mathrm{as}\quad t\to\infty.
\]
\item If $u(t,x)=e^{-at}$, $a>0$, then the asymptotic of $u^{E}(t,x)$
and its Cesaro limit are equal to 
\[
Ct^{-\alpha}\quad\mathrm{as}\;t\to\infty.
\]
\end{enumerate}
\end{thm}

The proof of Theorem \ref{thm:ISS} is essentially the contents of
Section \ref{sec:ISS} while the next theorem is shown in Section
\ref{sec:General-Classes}.
\begin{thm}
\label{thm:general-class}Let $u^{E}(t,x)$ be the subordination by
the density $G_{t}(\tau)$ associated to the classes \ref{eq:C1},
\ref{eq:C2}, and \ref{eq:C3} and $M_{t}(u^{E}(\cdot,x)):=\frac{1}{t}\int_{0}^{t}u^{E}(s,x)\,\mathrm{d}s$
the Cesaro mean of $u^{E}(t,x)$. 

\begin{enumerate}
\item Assume that $u(t,x)=t^{n}$, $n\ge0$. Then the asymptotic of the
Cesaro mean for the three classes are: 

\begin{description}
\item [{(C1).}] $M_{t}(u^{E}(\cdot,x))\sim Ct^{\alpha n}$ as $t\to\infty$, 
\item [{(C2).}] $M_{t}(u^{E}(\cdot,x))\sim C\log(t)^{n}$ as $t\to\infty$, 
\item [{(C3).}] $M_{t}(u^{E}(\cdot,x))\sim C\log(t)^{(1+s)n}$ as $t\to\infty$. 
\end{description}
\item If $u(t,x)=e^{-at}$, $a>0$, then the asymptotic of $M_{t}(u^{E}(\cdot,x))$
for the different classes are: 

\begin{description}
\item [{(C1).}] $M_{t}(u^{E}(\cdot,x))\sim Ct^{-\alpha}$ as $t\to\infty$, 
\item [{(C2).}] $M_{t}(u^{E}(\cdot,x))\sim C\log(t)^{-1}$ as $t\to\infty$, 
\item [{(C3).}] $M_{t}(u^{E}(\cdot,x))\sim C\log(t)^{-1-s}$ as $t\to\infty$. 
\end{description}
\end{enumerate}
\end{thm}

\section{Inverse Stable Subordinators}

\label{sec:ISS}In this section we consider two elementary solutions
of dynamical systems, namely $u(t)=u(t,x)=t^{n}$, $n\ge0$ and $u(t)=u(t,x)=e^{-at}$,
$a>0$, and investigate their subordination by the density $G_{t}(\tau)$
of inverse stable subordinator.

Define the function $u^{E}(t)=u^{E}(t,x)$ as the subordination of
$u(t)$ (of the above type) by the kernel $G_{t}(\tau)$, that is,
\begin{equation}
u^{E}(t):=\int_{0}^{\infty}u(\tau)G_{t}(\tau)\,\mathrm{d}\tau,\quad t\ge0.\label{eq:subordination}
\end{equation}
Our goal is to investigate the asymptotic behavior of $u^{E}(t)$.
At first we compute explicitly the function $u^{E}(t)$ by solving
the integral \eqref{eq:subordination} and obtain the time asymptotic.
Second we derive the Cesaro limit of $u^{E}(t)$, more precisely,
the asymptotic behavior of the Cesaro mean of $u^{E}(t)$ defined
by 
\begin{equation}
M_{t}(u^{E}(\cdot)):=\frac{1}{t}\int_{0}^{t}u^{E}(s)\,\mathrm{d}s.\label{Cesaro-mean}
\end{equation}
It turns out that both asymptotic behaviors for the two functions
$u(t)$ given above coincide. Therefore, for the random time change
associated to the inverse stable subordinator $E(t)$, $t\ge0$, the
asymptotic behavior of $u^{E}(t)$ is the same as the Cesaro limit.
On the other hand, using the Cesaro limit we may investigate a broad
class of subordinators. In Section \ref{sec:General-Classes} we investigate
the Cesaro limit for the classes \ref{eq:C1}, \ref{eq:C2}, and \ref{eq:C3}
while in this section concentrate in the spacial case of inverse stable
subordinators.

\subsection{Subordination of Monomials}

Let us consider at first the subordination of the function $u(t)=t^{n}$,
$n\ge0$. Hence, $u^{E}(t)$ is given by 
\begin{equation}
u^{E}(t)=\int_{0}^{\infty}\tau^{n}G_{t}(\tau)\,\mathrm{d}\tau.\label{eq:subordination-tn}
\end{equation}

It follows from \eqref{eq:LT-ISS-density} that $u^{E}(t)$ is explicitly
evaluated as 
\[
u^{E}(t)=(-1)^{n}\frac{\mathrm{d}^{n}}{\mathrm{d}\lambda^{n}}E_{\alpha}(-\lambda t^{\alpha})\big|_{\lambda=0}=\frac{n!}{\Gamma(\alpha n+1)}t^{\alpha n}.
\]
The last equality follows easily from the power series of the Mittag-Leffler
function 
\[
E_{\alpha}(z)=\sum_{n=1}^{\infty}\frac{z^{n}}{\Gamma(\alpha n+1)}.
\]
In addition, the asymptotic of the Mittag-Leffler function $E_{\alpha}$
that gives 
\begin{equation}
u^{E}(t)\sim Ct^{n\alpha}\quad\mathrm{as}\quad t\to\infty.\label{eq:asymp-direct}
\end{equation}

Now we turn to compute the asymptotic behavior of the Cesaro mean
of $u^{E}(t)$ with the help of the FKT theorem. To this end we define
the monotone function 
\begin{equation}
v(t):=\int_{0}^{t}u^{E}(s)\,\mathrm{d}s.\label{eq:monotone-mon}
\end{equation}
The Laplace-Stieltjes transform $w(\lambda)$ of $v(t)$ is given
by 
\[
w(\lambda):=\int_{0}^{\infty}e^{-\lambda t}\mathrm{d}v(t)=\int_{0}^{\infty}e^{-\lambda t}u^{E}(t)\,\mathrm{d}t=\int_{0}^{\infty}e^{-\lambda t}\int_{0}^{\infty}\tau^{n}G_{t}(\tau)\,\mathrm{d}\tau\,\mathrm{d}t.
\]
Using Fubini's theorem and equation \eqref{eq:LT-G-t} we obtain 
\[
w(\lambda)=\int_{0}^{\infty}\tau^{n}\int_{0}^{\infty}e^{-\lambda t}G_{t}(\tau)\,\mathrm{d}t\,\mathrm{d}\tau=\mathcal{K}(\lambda)\int_{0}^{\infty}\tau^{n}e^{-\tau\lambda\mathcal{K}(\lambda)}\,\mathrm{d}\tau.
\]
The r.h.s.~integral can be evaluated as 
\[
\int_{0}^{\infty}\tau^{n}e^{-\tau\lambda\mathcal{K}(\lambda)}\,\mathrm{d}\tau=(\lambda\mathcal{K}(\lambda))^{-(1+n)}n!
\]
which yields 
\begin{equation}
w(\lambda)=n!\lambda^{-(1+n)}\mathcal{K}(\lambda)^{-n}.\label{eq:LST-tn}
\end{equation}
On the other hand, for the stable subordinator we have $\mathcal{K}(\lambda)=\lambda^{\alpha-1}$,
cf. Example~\ref{exa:alpha-stable1}. Thus, we obtain 
\[
w(\lambda)=n!\lambda^{-(1+\alpha n)}=\lambda^{-\rho}L\left(\frac{1}{\lambda}\right),
\]
where $\rho=1+\alpha n$ and $L(x)=n!$ is a trivial slowly varying
function. Then Theorem \ref{thm:FKT-LST} yields 
\[
v(t)\sim Ct^{1+n\alpha}\quad\mathrm{as}\quad t\to\infty
\]
and this implies the following asymptotic behavior for the Cesaro
mean of $u^{E}(t)$ 
\begin{equation}
M_{t}(u^{E}(\cdot))=\frac{1}{t}\int_{0}^{t}u^{E}(s)\,\mathrm{d}s\sim Ct^{\alpha n}\quad\mathrm{as}\quad t\to\infty.\label{eq:CM-tn-ISS}
\end{equation}

\begin{rem}
In conclusion, we find that the asymptotic behavior of the subordination
$u^{E}(t)$ of any monomial by the density $G_{t}(\tau)$ (of the
inverse stable subordinator) as well as its Cesaro limit coincides.
Note also the slower decay of the subordination $u^{E}(t)$ compared
to $u(t)$ due to $0<\alpha<1$. 
\end{rem}

\subsection{Subordination of Decaying Exponentials }

Now we consider the solution $u(t)=e^{-at}$, $a>0$ and proceed to
study the asymptotic behavior of its subordination $u^{E}(t)$ by
the kernel $G_{t}(\tau)$. Again a direct computation is possible
in that case as well as the Cesaro mean.

Hence, the subordination $u^{E}(t)$ is given by 
\begin{equation}
u^{E}(t)=\int_{0}^{\infty}u(\tau)G_{t}(\tau)\,\mathrm{d}\tau=\int_{0}^{\infty}e^{-a\tau}G_{t}(\tau)\,\mathrm{d}\tau.\label{eq:subordination-exp}
\end{equation}
It follows from equation \eqref{eq:LT-ISS-density} that 
\begin{equation}
u^{E}(t)=E_{\alpha}(-at^{\alpha})\sim Ct^{-\alpha}\quad\mathrm{as}\quad t\to\infty.\label{eq:asymp-direct-exp}
\end{equation}

On the other hand, to derive the asymptotic behavior for the Cesaro
mean of $u^{E}(t)$ (with the help of Theorem~\ref{thm:FKT-LST})
we define the monotone function 
\begin{equation}
v(t):=\int_{0}^{t}u^{E}(s)\,\mathrm{d}s.\label{eq:monotone-exp}
\end{equation}
The Laplace-Stieltjes transform $w(\lambda)$ of $v(t)$ is equal
to 
\[
w(\lambda):=\int_{0}^{\infty}e^{-\lambda t}\mathrm{d}v(t)=\int_{0}^{\infty}e^{-\lambda t}u^{E}(t)\,\mathrm{d}t=\int_{0}^{\infty}e^{-\lambda t}\int_{0}^{\infty}e^{-a\tau}G_{t}(\tau)\,\mathrm{d}\tau\,\mathrm{d}t
\]
and using Fubini's theorem and equation \eqref{eq:LT-ISS-density}
we obtain 
\begin{equation}
w(\lambda)=\mathcal{K}(\lambda)\int_{0}^{\infty}e^{-\tau(a+\lambda\mathcal{K}(\lambda))}\,\mathrm{d}\tau=\frac{\mathcal{K}(\lambda)}{a+\lambda\mathcal{K}(\lambda)}.\label{eq:LST-exp}
\end{equation}
As $\mathcal{K}(\lambda)=\lambda^{\alpha-1}$ for the class \ref{eq:C1}
we may write $\tilde{v}(\lambda)$ as 
\[
w(\lambda)=\lambda^{-(1-\alpha)}\frac{1}{a+\lambda^{\alpha}}=\lambda^{-\rho}L\left(\frac{1}{\lambda}\right),\qquad\rho=1-\alpha,\quad L(t):=\frac{1}{a+t^{-\alpha}}.
\]
It is simple to verify that $L$ is a slowly varying function so that
we may use the FKT theorem to obtain 
\[
v(t)\sim Ct^{1-\alpha}\frac{1}{a+t^{-\alpha}}\quad\mathrm{as}\quad t\to\infty.
\]
Dividing both sides by $t$ leads to the asymptotic behavior of the
Cesaro mean of $u^{E}(t)$, that is, 
\begin{equation}
M_{t}(u^{E}(\cdot))=\frac{1}{t}\int_{0}^{t}u^{E}(s,x)\,\mathrm{d}s\sim C\frac{t^{-\alpha}}{a+t^{-\alpha}}\sim Ct^{-\alpha}\quad\mathrm{as}\;t\to\infty.\label{eq:CM-exp-ISS}
\end{equation}

\begin{rem}
We conclude that the asymptotic behavior $u^{E}(t)$ given in \eqref{eq:asymp-direct-exp}
coincides with the Cesaro limit of $u^{E}(t,x)$. In addition, we
notice that the starting function $u(t)=e^{-at}$ has an exponential
decay and its subordination has a slower decay, namely polynomial
decay.
\end{rem}

\section{Cesaro Limit for General Classes of Subordinators}

\label{sec:General-Classes}In this section we study the asymptotic
behavior of the subordination by the density $G_{t}(\tau)$ associated
to the classes \ref{eq:C1}, \ref{eq:C2}, and \ref{eq:C3}. Note
that Example \ref{exa:alpha-stable1} and \ref{exa:sum-two-stables}
belong to the class \ref{eq:C1}. As pointed out in Section \ref{sec:ISS}
here we only study the Cesaro limit of the subordination function
$u^{E}(t)$.

As in Section \ref{sec:ISS}, $u^{E}(t)$ is defined by 
\begin{equation}
u^{E}(t):=\int_{0}^{\infty}\tau^{n}G_{t}(\tau)\,\mathrm{d}\tau\label{eq:uE-tn}
\end{equation}
or 
\begin{equation}
u^{E}(t):=\int_{0}^{\infty}e^{-a\tau}G_{t}(\tau)\,\mathrm{d}\tau\label{eq:uE-exp}
\end{equation}
while $v(t)$ is defined by 
\[
v(t):=\int_{0}^{t}u^{E}(s)\,\mathrm{d}s.
\]
The density $G_{t}(\tau)$ in \eqref{eq:uE-tn} and \eqref{eq:uE-exp}
is associated to each class \ref{eq:C1}--\ref{eq:C3} described
above. We study the Cesaro limit of $u^{E}(t)$ for each class separately.

\subsection{Subordination by the Class (C1)}

At first we study the asymptotic behavior of $u^{E}(t)$ given by
\eqref{eq:uE-tn}. To this end we use equality \eqref{eq:LST-tn}
to obtain the Laplace-Stieltjes transform $w(\lambda)$ of the function
$v(t)$ as 
\[
w(\lambda):=\int_{0}^{\infty}e^{-\lambda t}\mathrm{d}v(t)=\lambda^{-(1+n)}(\mathcal{K}(\lambda))^{-n}n!.
\]
It follows from the behavior of $\mathcal{K}(\lambda)$ at $\lambda=0$
of the class \ref{eq:C1} that 
\[
w(\lambda)\sim\lambda^{-(1+\alpha n)}n!=\lambda^{-\rho}L\left(\frac{1}{\lambda}\right),
\]
where $\rho=1+\alpha n$ and $L(x)=n!$ is a slowly varying function.
It follows from the FKT theorem that 
\[
v(t)\sim Ct^{\rho}L(t)=Ct^{1+\alpha n}\quad\mathrm{as}\quad t\to\infty.
\]
This implies the Cesaro limit of $u^{E}(t)$ as 
\[
M_{t}(u^{E}(\cdot))\sim Ct^{\alpha n}\quad\mathrm{as}\quad t\to\infty.
\]
Note that this asymptotic is similar to the analogous for the inverse
stable subordinator, cf.~\eqref{eq:CM-tn-ISS}.

Let us now study the Cesaro limit of the function $u^{E}(t)$ given
in \eqref{eq:uE-exp}. Using the equality \eqref{eq:LST-exp} the
Laplace-Stieltjes transform $v(t)$ has the form 
\[
\tilde{v}(\lambda)=\frac{\mathcal{K}(\lambda)}{a+\lambda\mathcal{K}(\lambda)}.
\]
Replacing the local behavior of $\mathcal{K}(\lambda)$ at $\lambda=0$
for the class \ref{eq:C1} gives 
\[
\tilde{v}(\lambda)\sim\frac{\lambda^{\alpha-1}}{a+\lambda^{\alpha}}=\lambda^{-\rho}L\left(\frac{1}{\lambda}\right),
\]
where $\rho=1-\alpha$ and $L(x)=\frac{1}{1+ax^{-\alpha}}$. An applications
of the FKT theorem yields the asymptotic for $v(t)$, namely $v(t)\sim Ct^{\rho}L(t)$
as $t\to\infty$. Finally dividing both sides by $t$ gives the Cesaro
limit of $u^{E}(t)$, that is, 
\[
M_{t}(u^{E}(\cdot))\sim C\frac{t^{-\alpha}}{1+at^{-\alpha}}\sim Ct^{-\alpha}\quad\mathrm{as}\quad t\to\infty.
\]
Again, we obtain the same asymptotic as for the inverse stable subordinator,
see \eqref{eq:CM-exp-ISS}. In any case, since $0<\alpha<1$, the
time decaying is slower than the initial function $u(t)$.

\subsection{Subordination by the Class (C2)}

Assume that $u^{E}(t)$ is the subordination given in \eqref{eq:uE-tn}.
The Laplace-Stieltjes transform $w(\lambda)$ of $v(t)$ (cf.~equality
\eqref{eq:LST-tn}) has the form 
\[
w(\lambda):=\int_{0}^{\infty}e^{-\lambda t}\mathrm{d}v(t)=\lambda^{-(1+n)}(\mathcal{K}(\lambda))^{-n}n!.
\]
Using the behavior of $\mathcal{K}(\lambda)$ near $\lambda=0$ for
the class \ref{eq:C2} we obtain 
\[
w(\lambda)\sim\lambda^{-1}L\left(\frac{1}{\lambda}\right),
\]
where $L(x)=C\log(x)^{n}$, $C>0$, is a slowly varying function.
Then it follows from the FKT theorem that 
\[
v(t)\sim Ct\log(t)^{n}
\]
and as a result the asymptotic behavior for the Cesaro mean of $u^{E}(t)$
follows 
\[
M_{t}(u^{E}(\cdot))\sim C\log(t)^{n}\quad\mathrm{as}\quad t\to\infty.
\]

A similar analysis may be applied to study the asymptotic behavior
for the subordination $u^{E}(t)$ given in \eqref{eq:uE-exp}. The
Laplace-Stieltjes transform $w(\lambda)$ of the monotone function
$v(t)$ may be evaluated using equality \eqref{eq:LST-exp} to find
the following expression 
\[
w(\lambda)=\frac{\mathcal{K}(\lambda)}{a+\lambda\mathcal{K}(\lambda)}.
\]
Using the local behavior of $\mathcal{K}(\lambda)$ near $\lambda=0$
from class \ref{eq:C2} yields 
\[
w(\lambda)\sim\lambda^{-1}L\left(\frac{1}{\lambda}\right),
\]
where $L\left(x\right)=C\frac{\log\left(x\right)^{-1}}{a+C\log(x)^{-1}}$
which is a slowly varying function. Using the FKT theorem we obtain
the longtime behavior for the Cesaro mean of $u^{E}(t)$ as 
\[
M_{t}(u^{E}(\cdot))\sim C\frac{\log\left(t\right)^{-1}}{a+C\log(t)^{-1}}\sim C\log(t)^{-1}\quad\mathrm{as}\quad t\to\infty.
\]

\subsection{Subordination by the Class (C3)}

At first we study the subordination $u^{E}(t)$ given in \eqref{eq:uE-tn}
for the class \ref{eq:C3}. The Laplace-Stieltjes transform $w(\lambda)$
of the corresponding $v(t)$ is computed using equality \eqref{eq:LST-tn}
and we obtain 
\[
w(\lambda):=\int_{0}^{\infty}e^{-\lambda t}\mathrm{d}v(t)=\lambda^{-(1+n)}(\mathcal{K}(\lambda))^{-n}n!.
\]
Using the behavior of $\mathcal{K}(\lambda)$ near $\lambda=0$ for
the class \ref{eq:C3} yields 
\[
w(\lambda)\sim\lambda^{-1}L\left(\frac{1}{\lambda}\right),
\]
where $L(x)=C\log(x)^{(1+s)n}$, $C>0$, is a slowly varying function.
Then it follows from Theorem \ref{thm:FKT-LST} that 
\[
v(t)\sim Ct\log(t)^{(1+s)n}
\]
and dividing both sides by $t$ gives the asymptotic behavior for
the Cesaro mean of $u^{E}(t)$, namely 
\[
M_{t}(u^{E}(\cdot))\sim C\log(t)^{(1+s)n}\quad\mathrm{as}\quad t\to\infty.
\]

Let $u^{E}(t)$ be the subordination by $u(t)=e^{-at}$, $a>0$, that
is, equality \eqref{eq:uE-exp} with $G_{t}(\tau)$ from the class
\ref{eq:C3}. It follows from equality \eqref{eq:LST-exp} that the
Laplace-Stieltjes transform $w(\lambda)$ of $v(t)$ has the form
\[
w(\lambda)=\frac{\mathcal{K}(\lambda)}{a+\lambda\mathcal{K}(\lambda)}.
\]
Using the local behavior of $\mathcal{K}(\lambda)$ near $\lambda=0$
from class \ref{eq:C3} yields 
\[
w(\lambda)\sim\lambda^{-1}L\left(\frac{1}{\lambda}\right),\qquad L\left(x\right)=C\frac{\log\left(x\right)^{-1-s}}{a+C\log(x)^{-1-s}},
\]
where $C,s>0$. As the function $L$ is slowly varying at infinity,
then by the FKT theorem we obtain the asymptotic behavior for the
Cesaro mean of $u^{E}(t)$ as 
\[
M_{t}(u^{E}(\cdot))\sim C\frac{\log\left(t\right)^{-1-s}}{a+C\log(t)^{-1-s}}\sim C\log(t)^{-1-s}\quad\mathrm{as}\quad t\to\infty.
\]

\subsection*{Acknowledgments}

This work has been partially supported by Center for Research in Mathematics
and Applications (CIMA) related with the Statistics, Stochastic Processes
and Applications (SSPA) group, through the grant UIDB/MAT/04674/2020
of FCT-Funda{\c c\~a}o para a Ci{\^e}ncia e a Tecnologia, Portugal.

The financial support by the Ministry for Science and Education of
Ukraine through Project 0119U002583 is gratefully acknowledged.

\end{document}